\documentclass[12pt]{article}
\usepackage{amssymb}
\usepackage{amsmath}
\newcommand{\la}{\lambda}
\newcommand{\lap}{\mbox{$\bigtriangleup$}}
\newcommand{\grad}{\mbox{$\bigtriangledown$}}
\newcommand{\ra}{{\mbox{$\rightarrow$}}}
\newcommand{\be}{\begin{equation}}
\newcommand{\ee}{\end{equation}}
\newtheorem{mrem}{Remark}
\newtheorem{mthm}{Theorem}
\newtheorem{mpro}{Proposition}

\newtheorem{mcor}{Corollary}
\newtheorem{thm}{Theorem}[section]
\newtheorem{pro}{Proposition}[section]
\newtheorem{lem}{Lemma}[section]

\newtheorem{rem}{Remark}[section]
\begin{document}

\title{A direct method of moving planes for the fractional Laplacian}

\author{Wenxiong Chen \thanks{Partially supported by the Simons Foundation Collaboration Grant for Mathematicians 245486.} \quad  Congming Li \thanks{ Corresponding author, Partially supported by NSF DMS-1405175 and NSFC-11271166.} \quad Yan Li  }

\date{\today}
\maketitle
\begin{abstract} In this paper, we develop a direct {\em method of moving planes} for the fractional Laplacian. Instead of using the conventional extension method introduced by Caffarelli and Silvestre, we work directly on the non-local operator. Using the integral defining the fractional Laplacian, by an elementary approach, we first obtain the key ingredients needed in the {\em method of moving planes} either in a bounded domain or in the whole space, such as {\em strong maximum principles for anti-symmetric functions}, {\em narrow region principles}, and {\em decay at infinity}. Then, using simple examples, semi-linear equations involving the fractional Laplacian,  we illustrate how this new {\em method of moving planes} can be employed to obtain symmetry and non-existence of positive solutions.

We firmly believe that the ideas and methods introduced here can be conveniently applied to study a variety of nonlocal problems with more general operators and more general nonlinearities.
\end{abstract}
\bigskip

\bigskip

\section{Introduction}

The fractional Laplacian in $R^n$ is a nonlocal pseudo-differential operator, assuming the form
\begin{equation}
(-\Delta)^{\alpha/2} u(x) = C_{n,\alpha} \, \lim_{\epsilon \ra 0} \int_{\mathbb{R}^n\setminus B_{\epsilon}(x)} \frac{u(x)-u(z)}{|x-z|^{n+\alpha}} dz,
\label{op}
\end{equation}
 where $\alpha$ is any real number between $0$ and $2$.
This operator is well defined in $\cal{S}$, the Schwartz space of rapidly decreasing $C^{\infty}$
functions in $\mathbb{R}^n$. In this space, it can also be equivalently defined in terms of the Fourier transform
$$ \widehat{(-\Delta)^{\alpha/2} u} (\xi) = |\xi|^{\alpha} \hat{u} (\xi), $$
where $\hat{u}$ is the Fourier transform of $u$. One can extend this operator to a wider space of functions.

Let
$$L_{\alpha}=\{u: \mathbb{R}^n\rightarrow \mathbb{R} \mid \int_{\mathbb{R}^n}\frac{|u(x)|}{1+|x|^{n+\alpha}} \, d x <\infty\}.$$
Then it is easy to verify that for $u \in L_{\alpha}  \cap C_{loc}^{1,1}$, the integral on the right hand side of (\ref{op}) is well defined. Throughout this paper, we consider the fractional Laplacian in this setting.

The non-locality of the fractional Laplacian makes it difficult to investigate. To circumvent this difficulty,
Caffarelli and Silvestre \cite{CS} introduced the {\em extension method} that reduced this nonlocal problem into a local one in higher dimensions. For a function $u:\mathbb{R}^n \ra \mathbb{R}$, consider the extension $U:\mathbb{R}^n\times[0, \infty) \ra \mathbb{R}$ that satisfies
$$\left\{\begin{array}{ll}
div(y^{1-\alpha} \nabla U)=0, & (x,y) \in \mathbb{R}^n\times[0, \infty),\\
U(x, 0) = u(x).
\end{array}
\right.
$$
Then
$$(-\lap)^{\alpha/2}u (x) = - C_{n,\alpha} \displaystyle\lim_{y \ra 0^+}
y^{1-\alpha} \frac{\partial U}{\partial y},  \;\; x \in \mathbb{R}^n.$$

This {\em extension method } has been applied successfully to study equations involving the fractional Laplacian, and a series of fruitful results have been obtained (see \cite{BCPS} \cite{CZ}    and the references therein).

In \cite{BCPS}, among many interesting results, when the authors considered the properties of the positive solutions for
\begin{equation} (-\lap)^{\alpha/2} u = u^p (x), \;\; x \in \mathbb{R}^n,
\label{e}
\end{equation}
they first used the above {\em extension method} to reduce the nonlocal problem into a local one for $U(x,y)$ in one higher dimensional half space $\mathbb{R}^n\times[0, \infty)$, then applied the {\em method of moving planes} to show the symmetry of $U(x,y)$ in $x$, and hence derived the non-existence in the subcritical case:
\begin{mpro} (Brandle-Colorado-Pablo-Sanchez)
Let $1 \leq \alpha <2$. Then the problem
\begin{equation}
\left\{\begin{array}{ll}
div(y^{1-\alpha} \nabla U)=0, & (x,y) \in \mathbb{R}^n\times[0, \infty),\\
\displaystyle-\lim_{y \ra 0^+} y^{1-\alpha} \frac{\partial U}{\partial y} = U^p (x,0), & x \in \mathbb{R}^n
\end{array}
\right.
\label{U}
\end{equation}
has no positive bounded solution provided $p < (n+\alpha)/(n-\alpha).$
\end{mpro}

They then took trace to obtain
\begin{mcor} Assume that $1\leq \alpha < 2$ and $1<p < \frac{n+\alpha}{n-\alpha}$. Then equation (\ref{e}) possesses no bounded positive solution.
\label{mcor}
\end{mcor}

A similar {\em extension method} was adapted in \cite{CZ} to obtain the nonexistence of positive solutions for an indefinite fractional problem:

\begin{mpro} (Chen-Zhu)
Let $1 \leq \alpha <2$ and $1<p<\infty$. Then the equation
\begin{equation}
(-\lap)^{\alpha/2} u = x_1 u^p, \;\; x \in \mathbb{R}^n
\label{x1e}
\end{equation}
possesses no positive bounded solutions.
\end{mpro}

The common restriction $\alpha \geq 1$ is due to the approach that they need to carry the {\em method of moving planes} on the solutions $U$ of the extended problem
\begin{equation}
div(y^{1-\alpha} \nabla U)=0, \;\; (x,y) \in \mathbb{R}^n\times[0, \infty).
\label{eeq}
\end{equation}
 Due to technical restriction, they have to assume $\alpha \geq 1$. It seems that this condition cannot be weakened if one wants to carry the {\em method of moving planes} on extended equation (\ref{eeq}).

{\em Then what happens in the case $0< \alpha < 1$? }

Actually, this case can be treated by considering the corresponding integral equation. In \cite{CLO} \cite{CLO1}, the authors showed that if $u \in H^{\alpha/2}(\mathbb{R}^n)$ is a positive weak solution of (\ref{e}), then it also satisfies the integral equation
\begin{equation}
 u(x) = C \int_{\mathbb{R}^n} \frac{1}{|x-y|^{n-\alpha}} u^p(y) d y .
 \label{ie}
 \end{equation}
Applying the {\em method of moving planes in integral forms}, they obtained the radial symmetry in the critical case and non-existence in the subcritical case for positive solutions of (\ref{ie}).

Under the weaker condition that $u \in L_{\alpha}(\mathbb{R}^n)$, the equivalence between pseudo differential equation (\ref{e}) and integral equation was also established in \cite{ZCCY} by employing a Liouville theorem for $\alpha$-harmonic functions.

In the case of more general nonlinearity, for instance, when considering
\begin{equation}
(- \lap)^{\alpha/2} u = f(x,u), \;\; x \in \mathbb{R}^n ,
\label{ef}
\end{equation}
in order to show that a positive solution of (\ref{ef}) also solves
\begin{equation}
 u(x) = C \int_{\mathbb{R}^n} \frac{1}{|x-y|^{n-\alpha}} f(x,u(y)) d y ,
 \label{ief}
 \end{equation}
 so far one needs to assume that $f(x,u)$ is nonnegative, which is not satisfied by the right hand side of equation (\ref{x1e}). Hence in this situation, the integral equation approach renders powerless.

 Another technical restriction in carrying out the method of moving planes on the integral equation is that both $f(x, u)$ and $\frac{\partial f}{\partial u}$ must be monotone increasing in $u$, which may not be necessary if one directly works on pseudo differential equation (\ref{ef}).

In summary, either by {\em extension}  or by integral equations, one needs to impose extra conditions on the solutions, which would not be necessary if we consider the pseudo differential equation directly. Moreover, for equations involving fully nonlinear operators, such as 
$$F_{\alpha}(u) = f(x,u) $$
where 
\begin{equation} 
F_{\alpha}(u(x)) = C_{n,\alpha} \, \lim_{\epsilon \ra 0} \int_{\mathbb{R}^n\setminus B_{\epsilon}(x)} \frac{G(u(x)-u(z))}{|x-z|^{n+\alpha}} dz,
\label{F}
\end{equation}
with $G(\cdot)$ being a Lipschitz continuous function (see \cite{CS1}), so far as we know, there has neither  been any corresponding {\em extension methods} nor equivalent integral equations that one can work at.

{\em Then can one carry on the method of moving planes directly on nonlocal equations? }

The main objective of this paper is to answer this question affirmatively.

Actually, the first partial answer was provided in \cite{JW} by Jarohs and Weth. There they introduced {\em antisymmetric maximum principles} and applied them to carry on the {\em method of moving planes} directly on nonlocal problems to show the symmetry of solutions. The operators they considered are quite general, however, their maximum principles only apply to bounded regions $\Omega$, and they only considered weak solutions defined by $H^{\alpha/2}(\Omega)$ inner product.

In this paper, we will develop a systematical approach to carry on the {\em method of moving planes} for nonlocal problems, either on bounded or unbounded domains. For local elliptic operators, these kinds of
approaches were introduced decades ago in the first two authors' paper \cite{CL} and then summarized in their book \cite{CL1}, among which the {\em narrow region principle} and the {\em decay at infinity} have been applied extensively by many researchers to solve various problems. A parallel system for the fractional Laplacian will be established here by very elementary methods, so that it can be conveniently applied to various nonlocal problems. This will be accomplished in Section 2. The main theorems and how they fit in the framework of the method of moving planes are illustrated in the following.
\bigskip

\leftline{\bf Key Ingredients in the Method of Moving Planes}
\smallskip

As usual, let
$$T_{\lambda} =\{x \in \mathbb{R}^{n}|\; x_1=\lambda, \mbox{ for some } \lambda\in \mathbb{R}\}$$
be the moving planes,
$$\Sigma_{\lambda} =\{x \in \mathbb{R}^{n} | \, x_1<\lambda\}$$
be the region to the left of the plane, and
$$ x^{\lambda} =(2\lambda-x_1, x_2, ..., x_n)$$
be the reflection of $x$ about the plane $T_{\lambda}$.

Assume that $u$ is a solution of pseudo differential equation (\ref{e}) or (\ref{ef}). To compare the values of $u(x)$ with $u(x^{\lambda})$, we denote
$$w_{\lambda} (x) = u(x^{\lambda}) - u(x) .$$

The first step is to show that for $\lambda$ sufficiently negative, we have
\begin{equation}
w_{\lambda}(x) \geq 0 , \;\; x \in \Sigma_{\lambda} .
\label{w}
\end{equation}
This provides a starting point to move the plane. Then in the second step, we move the plane to the right as long as inequality (\ref{w}) holds to its limiting position to show that $u$ is symmetric
about the limiting plane. A maximum principle is used to prove (\ref{w}). Since $w_{\lambda}$ is an anti-symmetric function:
$$w_{\lambda}(x) = - w_{\lambda}(x^{\lambda}),$$
we first prove (for simplicity of notation, in the following, we denote $w_{\lambda}$ by $w$ and $\Sigma_{\lambda}$ by $\Sigma$. )

\begin{mthm} ( Maximum Principle for Anti-symmetric Functions.)

Let $\Omega$ be a bounded domain in $\Sigma$.
Assume that $w\in L_\alpha\cap C_{loc}^{1,1}(\Omega)$ and is lower semi-continuous on $\bar{\Omega}$. If
$$
\left\{\begin{array}{ll}
(-\lap)^{\alpha/2}w(x) \geq 0  &\mbox{ in } \Omega,\\
w(x) \geq0&\mbox{ in }  \Sigma \backslash\Omega,\\
w(x^{\lambda})=-w(x)  &\mbox{ in } \Sigma,
\end{array}
\right.
$$
then
$$
w(x) \geq0 \mbox{ in } \Omega.
$$

Furthermore, if $w = 0$ at some point in $\Omega$, then
$$ w(x) = 0 \; \mbox{ almost everywhere in }  \mathbb{R}^n. $$

These conclusions hold for unbounded region $\Omega$ if we further assume that
$$\underset{|x| \ra \infty}{\underline{\lim}} w(x) \geq 0 .$$

\end{mthm}

In many cases, $w$ may not satisfy the equation
$$ (-\lap)^{\alpha/2} w \geq 0 $$
as required in the previous theorem. However one can derive that
$$ (-\lap)^{\alpha/2} w  + c(x) w (x) \geq 0 $$
for some function $c(x)$ depending on $u$. If $c(x)$ is nonnegative, it is easy to see that the {\em maximum principle} is still valid; however this is not the case in practice.
Fortunately, in the process of moving planes, each time we only need to move $T_{\lambda}$ a little bit to the right, hence the increment of $\Sigma_{\lambda}$ is a narrow region, and a {\em maximum principle} is easier to hold in a narrow region provided $c(x)$ is not ``too negative'', as you will see below.

\begin{mthm} ( Narrow Region Principle.)

Let $\Omega$ be a bounded narrow region in $\Sigma$, such that it is contained in
 $$\{x | \; \lambda-\delta<x_1<\lambda \, \}$$
 with small $\delta$. Suppose  that $w\in L_\alpha\cap C_{loc}^{1,1}(\Omega)$ and is lower semi-continuous on $\bar{\Omega}$. If
 $c(x)$ is bounded from below in $\Omega$  and
$$
\left\{\begin{array}{ll}
(-\lap)^{\alpha/2}w(x) +c(x)w(x)\geq0  &\mbox{ in } \Omega,\\
w(x) \geq0&\mbox{ in }  \Sigma \backslash\Omega,\\
w(x^{\lambda})=-w(x)  &\mbox{ in } \Sigma,
\end{array}
\right.
$$

then for sufficiently small $\delta$, we have
$$
w(x) \geq0 \mbox{ in } \Omega.
$$

Furthermore, if $w = 0$ at some point in $\Omega$, then
$$ w(x) = 0 \; \mbox{ almost everywhere in }  \mathbb{R}^n. $$

These conclusions hold for unbounded region $\Omega$ if we further assume that
$$\underset{|x| \ra \infty}{\underline{\lim}} w(x) \geq 0 .$$
\end{mthm}

As one will see from the proof of this theorem, the contradiction arguments are conducted at a negative minimum of $w$. Hence when working on an unbounded domain, one needs to rule out the possibility that such minima would ``leak'' to infinity. This can be done when $c(x)$ decays
``faster'' than $1/|x|^{\alpha}$ near infinity.

\begin{mthm} ( Decay at Infinity.)

Let $\Omega$ be an unbounded region in $\Sigma$.
Assume $w \in L_{\alpha} \cap C_{loc}^{1,1}(\Omega)$ is a solution of
$$
\left\{\begin{array}{ll}
(-\lap)^{\alpha/2}w(x) +c(x)w(x)\geq0  &\mbox{ in } \Omega,\\
w(x) \geq0&\mbox{ in } \Sigma \backslash\Omega,\\
w(x^{\lambda})=-w(x)  &\mbox{ in } \Sigma,
\end{array}
\right.
$$

with
$$
\underset{|x|\rightarrow \infty}{\underline{\lim}}|x|^{\alpha}c(x)\geq0,
$$
then there exists a constant $R_0>0$ ( depending on $c(x)$, but independent of $w$ ), such that if
$$
w(x^0)=\underset{\Omega}{\min}\;w(x)<0,
$$
 then
 $$
|x^0|\leq R_0.
$$
\label{mthmdi}
\end{mthm}

\leftline{\bf Applications of the Method of Moving Planes--Examples}
\smallskip

In Section 3, we will use several examples to illustrate how the key ingredients obtained in Section 2 can be used in the {\em method of moving planes} to establish symmetry and monotonicity of positive solutions.

We first consider
\begin{equation} (-\lap)^{\alpha/2} u = u^p (x), \;\; x \in \mathbb{R}^n,
\label{e1}
\end{equation}
and prove
\begin{mthm}
Assume that $0 < \alpha < 2$ and $u \in L_{\alpha} \cap C^{1,1}_{loc}$ is a nonnegative solution of
equation (\ref{e1}). Then

(i) In the critical case $p = \frac{n+\alpha}{n-\alpha}$, $u$ is radially symmetric and monotone decreasing about some point.

(ii) In the subcritical case $1<p< \frac{n+\alpha}{n-\alpha}$, $u \equiv 0$.
\label{mthm1}
\end{mthm}

\begin{mrem}
As compared to Corollary \ref{mcor}, we relaxed the condition $1\leq \alpha < 2$ to $0<\alpha<2$, and we dropped the global boundedness assumption on $u$. The local $C^{1,1}$ property for the solutions of (\ref{e1}) can be obtained by a standard regularity argument if $u \in L_{\alpha}$ is a solution in the sense of distribution.
\end{mrem}
\medskip

Then we investigate the same equation on upper half space with the Dirichlet condition:
\begin{equation}
\left\{\begin{array}{ll}
(-\lap)^{\alpha/2} u = u^p (x), & x \in \mathbb{R}_+^n, \\
u(x) \equiv 0 , & x \not{\in} \mathbb{R}_+^n .
\end{array}
\right.
\label{e1h}
\end{equation}

\begin{mthm}
Assume that $0 < \alpha < 2$ and $u \in L_{\alpha} \cap C^{1,1}_{loc}$ is a nonnegative solution of
equation (\ref{e1h}). Then
in the subcritical and critical case $1<p \leq \frac{n+\alpha}{n-\alpha}$, we have $u \equiv 0$.
\label{mthm1h}
\end{mthm}

 Next we study positive solutions for the nonlinear Schr$\ddot{o}$dinger equation with fractional diffusion
\begin{equation}
(-\lap)^{\alpha/2} u + u = u^p , \;\; x \in R^n .
\label{se}
\end{equation}
We prove

\begin{mthm}
Assume that $u \in L_{\alpha}\cap C_{loc}^{1,1}$ is a positive solution of (\ref{se}) with $1<p<\infty$. If
$$
\lim_{|x| \ra \infty} u(x) = a < \left( \frac{1}{p} \right)^{\frac{1}{p-1}} ,
$$
then $u$ must be radially symmetric and monotone decreasing about some point in $\mathbb{R}^n$.
\end{mthm}

\begin{mrem}
In \cite{FLe} and \cite{FLS}, Frank, Lenzmann, and Silvestre obtained the radial symmetry and uniqueness of the solution for equation (\ref{se}) without positivity assumption on $u$. However, they required that $u \in H^{\alpha/2}(\mathbb{R}^n)$ and $p < \frac{n+\alpha}{n-\alpha}$.
\end{mrem}

Finally, as a byproduct, we consider the problem with more general nonlinearity on a bounded domain:
\begin{equation}
\left\{ \begin{array}{ll}
(-\lap)^{\alpha/2} u(x)  = f(u(x)) , & x \in B_1(0), \\
u(x) = 0 , & x \not{\in} B_1(0).
\end{array}
\right.
\label{fu}
\end{equation}
We have
\begin{mthm}
Assume that $u \in L_{\alpha}\cap C_{loc}^{1,1}(B_1(0))$ is a positive solution of (\ref{fu}) with $f(\cdot)$ being Lipschitz continuous. Then $u$ must be radially symmetric and monotone decreasing about the origin.
\label{mthmfu}
\end{mthm}

Recently, the first author and G. Li \cite{CLg} applied the {\em method of moving planes} introduced in this paper to study the nonlinear problem involving fully nonlinear operator 
$$\left\{\begin{array}{ll} 
F_{\alpha}(u(x)) = f(x,u) & x \in \Omega \\
u \equiv 0 & x \not{\in} \Omega ,
\end{array}
\right.
$$
where $F_{\alpha}(\cdot)$ is defined in (\ref{F}). Radial symmetry and non-existence of solutions are established when $\Omega$ is a unit ball, or a half space, or the whole space. It is interesting to point out that, as $\alpha \ra 2$, 
$$ F_{\alpha}(u(x)) \ra -\lap u(x) + a |\grad u(x)|^2. $$

\begin{mrem}
(i) Theorem \ref{mthmfu} is an extension of the elegant result of Gidas, Ni, and Nirenberg \cite{GNN} on the Laplacian to the fractional Laplacian.

(ii) A similar result has been obtained in \cite{JW} under slightly different regularity assumptions on $u$.

\end{mrem}

{\bf Note.} After this paper has been posted on the arXiv, Xiong brought to our attention that they proved a maximum principle for odd solutions of a nonlocal parabolic equation in \cite{JLX}, and in their other paper \cite{JX} (Theorem 1.8), they obtained classifications of solutions for (\ref{U}) in the critical case for $0 < \alpha < 2.$ From \cite{JX}, we also learnt that in \cite{FLe1}, Frank and Lenzmann obtained a strong maximum principle and a Hopf lemma for odd solutions to nonlocal elliptic equations.

For more articles concerning the method of moving planes for nonlocal equations, mainly for integral equations, please see \cite{FL} \cite{Ha} \cite{HLZ} \cite{HWY} \cite{Lei} \cite{LLM} \cite{LZ} \cite{LZ1} \cite{LZ2} \cite{MC} \cite{MZ} and the references therein.

\section{Various Maximum Principles}

\subsection{A Maximum Principle for Anti-symmetric Functions}

We first provide a simpler proof for a well-known maximum principle for $\alpha$-super harmonic functions.

\begin{thm} (Maximum Principle)
Let $\Omega$ be a bounded domain in $\mathbb{R}^n$. Assume that $u\in L_\alpha\cap C_{loc}^{1,1}(\Omega)$ and is lower semi-continuous on $\bar{\Omega}$. If
\be
\left\{\begin{array}{ll}
(-\lap)^{\alpha/2}u(x) \geq0  &\mbox{ in } \Omega,\\
u(x) \geq0&\mbox{ in } \mathbb{R}^{n}\backslash\Omega,
\end{array}
\right.
\label{1}
\ee
then
\be
u(x) \geq0 \mbox{ in } \Omega.\label{2}
\ee

If $u = 0$ at some point in $\Omega$, then
$$ u(x) = 0 \; \mbox{ almost everywhere in }  \mathbb{R}^n. $$

These conclusions hold for unbounded region $\Omega$ if we further assume that
$$\underset{|x| \ra \infty}{\underline{\lim}} u(x) \geq 0 .$$
\label{thm2.1}
\end{thm}

\begin{rem}
This {\em maximum principle} has been established by Silvestre in \cite{Si} without the condition
$u \in C_{loc}^{1,1}(\Omega)$. Here we provide a much more elementary and simpler proof.
\end{rem}

$\mathbf{Proof.}$ \quad If (\ref{2}) does not hold,
then the lower semi-continuity of $u$ on $\bar{\Omega}$ indicates that there exists a
$x^0\in\bar{\Omega}$ such that
$$u(x^0)=\underset{\bar{\Omega}}{\min} \,u<0.$$
And one can further deduce from condition (\ref{1}) that $x^0$ is in the interior of $\Omega$.

Then it follows that
\begin{eqnarray*}
(-\lap)^{\alpha/2}u(x^0)&=&C_{n,\alpha}PV\int_{\mathbb{R}^n}\frac{u(x^0)-u(y)}
{|x^0-y|^{n+\alpha}}dy\\
&\leq&C_{n,\alpha}\int_{\mathbb{R}^n\backslash\Omega}\frac{u(x^0)-u(y)}
{|x^0-y|^{n+\alpha}}dy\\
&<&0,
\end{eqnarray*}
which contradicts inequality (\ref{1}). This verifies (\ref{2}).

If at some point $x^o \in \Omega$, $u(x^o) = 0$, then from
$$ 0 \leq (-\lap)^{\alpha/2}u(x^o) = C_{n,\alpha}PV\int_{\mathbb{R}^n}\frac{-u(y)}
{|x^o-y|^{n+\alpha}}dy $$
and $u \geq 0$, we must have
$$ u(x) = 0 \; \mbox{ almost everywhere in }  \mathbb{R}^n. $$

This completes the proof.
\medskip

Then we introduce a maximum principle for anti-symmetric functions.

\begin{thm}
Let $T$ be a hyperplane in $\mathbb{R}^{n}$. Without loss of generality,
we may assume that
$$T=\{x \in \mathbb{R}^{n}|\; x_1=\lambda, \mbox{ for some } \lambda\in \mathbb{R}\}.$$
Let $$\tilde{x}=(2\lambda-x_1, x_2, ..., x_n)$$
be the reflection of $x$ about the plane $T$. Denote
$$H=\{x \in \mathbb{R}^{n} | \, x_1<\lambda\} \; \mbox{ and } \; \tilde{H}=\{x | \, \tilde{x} \in H\},$$
Let $\Omega$ be a bounded domain in $H$.
Assume that $u\in L_\alpha\cap C_{loc}^{1,1}(\Omega)$ and is lower semi-continuous on $\bar{\Omega}$. If
\be
\left\{\begin{array}{ll}
(-\lap)^{\alpha/2}u(x) \geq0  &\mbox{ in } \Omega,\\
u(x) \geq0&\mbox{ in } H\backslash\Omega,\\
u(\tilde{x})=-u(x)  &\mbox{ in } H,
\end{array}
\right.
\label{3}
\ee
then
\be
u(x) \geq0 \mbox{ in } \Omega.\label{4}
\ee

Furthermore, if $u = 0$ at some point in $\Omega$, then
$$ u(x) = 0 \; \mbox{ almost everywhere in }  \mathbb{R}^n. $$

These conclusions hold for unbounded region $\Omega$ if we further assume that
$$\underset{|x| \ra \infty}{\underline{\lim}} u(x) \geq 0 .$$

\label{thmAS}
\end{thm}

$\mathbf{Proof.}$ \quad If (\ref{4}) does not hold,
then the lower semi-continuity of $u$ on $\bar{\Omega}$ indicates that there exists a
$x^0\in\bar{\Omega}$ such that
$$u(x^0)=\underset{\bar{\Omega}}{\min} \,u<0.$$
And one can further deduce from condition (\ref{3}) that $x^0$ is in the interior of $\Omega$.

It follows that
\begin{eqnarray*}
(-\lap)^{\alpha/2}u(x^0)&=&C_{n,\alpha}PV\int_{\mathbb{R}^n}\frac{u(x^0)-u(y)}
{|x^0-y|^{n+\alpha}}dy\\
&=&C_{n,\alpha}PV\left\{\int_H \frac{u(x^0)-u(y)}
{|x^0-y|^{n+\alpha}}dy+\int_{\tilde{H}} \frac{u(x^0)-u(y)}
{|x^0-y|^{n+\alpha}}dy \right\}\\
&=&C_{n,\alpha}PV\left\{\int_H \frac{u(x^0)-u(y)}
{|x^0-y|^{n+\alpha}}dy+\int_H \frac{u(x^0)-u(\tilde{y})}
{|x^0-\tilde{y}|^{n+\alpha}}dy\right\}\\
&=&C_{n,\alpha}PV\left\{\int_H \frac{u(x^0)-u(y)}
{|x^0-y|^{n+\alpha}}dy+\int_H \frac{u(x^0)+u(y)}
{|x^0-\tilde{y}|^{n+\alpha}}dy\right\}\\
&\leq&
C_{n,\alpha}\int_H
\left\{\frac{u(x^0)-u(y)}{|x^0-\tilde{y}|^{n+\alpha}}
+\frac{u(x^0)+u(y)}{|x^0-\tilde{y}|^{n+\alpha}}\right\}dy\\
&=&C_{n,\alpha}\int_H \frac{2u(x^0)}
{|x^0-\tilde{y}|^{n+\alpha}}dy\\
&<&0,
\end{eqnarray*}
which contradicts inequality (\ref{3}). This verifies (\ref{4}).

Now we have shown that $u \geq 0$ in $\mathbb{R}^n$. If there is some point $x^o \in \Omega$, such that $u(x^o) =0$, then from
$$ 0 \leq (-\lap)^{\alpha/2}u(x^o)= C_{n,\alpha}PV \int_{\mathbb{R}^n}\frac{-u(y)}
{|x^o-y|^{n+\alpha}}dy,$$
we derive immediately that
$$ u(x) = 0 \; \mbox{ almost everywhere in }  \mathbb{R}^n. $$
This completes the proof.

\subsection{Narrow Region Principle}

\begin{thm}
Let $T$ be a hyperplane in $\mathbb{R}^{n}$. Without loss of generality,
we may assume that
$$T=\{x=(x_1,x') \in \mathbb{R}^{n} | \; x_1=\lambda,\; \mbox{ for some } \lambda\in \mathbb{R}\}.$$
Let $$\tilde{x}=(2\lambda-x_1, x_2, ..., x_n),$$ $$H=\{x \in \mathbb{R}^{n} | \; x_1<\lambda\},\quad \tilde{H}=\{x | \; \tilde{x} \in H\}.$$
 Let $\Omega$ be a bounded narrow region in $H$, such that it is contained in $\{x | \; \lambda-l<x_1<\lambda \, \}$ with small $l$. Suppose  that $u\in L_\alpha\cap C_{loc}^{1,1}(\Omega)$ and is lower semi-continuous on $\bar{\Omega}$. If
 $c(x)$ is bounded from below in $\Omega$  and
\be
\left\{\begin{array}{ll}
(-\lap)^{\alpha/2}u(x) +c(x)u(x)\geq0  &\mbox{ in } \Omega,\\
u(x) \geq0&\mbox{ in } H\backslash\Omega,\\
u(\tilde{x})=-u(x)  &\mbox{ in } H,
\end{array}
\right.
\label{5}
\ee
then for sufficiently small $l$, we have
\be
u(x) \geq0 \mbox{ in } \Omega.\label{6}
\ee

Furthermore, if $u = 0$ at some point in $\Omega$, then
$$ u(x) = 0 \; \mbox{ almost everywhere in }  \mathbb{R}^n. $$

These conclusions hold for unbounded region $\Omega$ if we further assume that
$$\underset{|x| \ra \infty}{\underline{\lim}} u(x) \geq 0 .$$
\label{thmNR}
\end{thm}

$\mathbf{Proof.}$ \quad If (\ref{6}) does not hold,
then the lower semi-continuity of $u$ on $\bar{\Omega}$ indicates that there exists an
$x^0\in\bar{\Omega}$ such that
$$u(x^0)=\underset{\bar{\Omega}}{\min} \,u<0.$$
And one can further deduce from condition (\ref{5}) that $x^0$ is in the interior of $\Omega$.

Then it follows that
\begin{eqnarray*}
(-\lap)^{\alpha/2}u(x^0)&=&C_{n,\alpha}PV\int_{\mathbb{R}^n}\frac{u(x^0)-u(y)}
{|x^0-y|^{n+\alpha}}dy\\
&=&C_{n,\alpha}PV\left\{\int_H \frac{u(x^0)-u(y)}
{|x^0-y|^{n+\alpha}}dy+\int_{\tilde{H}} \frac{u(x^0)-u(y)}
{|x^0-y|^{n+\alpha}}dy \right\}\\
&=&C_{n,\alpha}PV\left\{\int_H \frac{u(x^0)-u(y)}
{|x^0-y|^{n+\alpha}}dy+\int_H \frac{u(x^0)-u(\tilde{y})}
{|x^0-\tilde{y}|^{n+\alpha}}dy\right\}\\
&=&C_{n,\alpha}PV\left\{\int_H \frac{u(x^0)-u(y)}
{|x^0-y|^{n+\alpha}}dy+\int_H \frac{u(x^0)+u(y)}
{|x^0-\tilde{y}|^{n+\alpha}}dy\right\}\\
&\leq&
C_{n,\alpha}\int_H
\left\{\frac{u(x^0)-u(y)}{|x^0-\tilde{y}|^{n+\alpha}}
+\frac{u(x^0)+u(y)}{|x^0-\tilde{y}|^{n+\alpha}}\right\}dy\\
&=&C_{n,\alpha}\int_H \frac{2u(x^0)}
{|x^0-\tilde{y}|^{n+\alpha}}dy.
\end{eqnarray*}

Let $D=\{y|l<y_1- x^0_1<1,\; |y'-(x^0)'|<1\}$, $s=y_1- x^0_1$,
$\tau=|y'-(x^0)'|$ and $\omega_{n-2}=|B_1(0)|$ in $R^{n-2}$. Now we have
\begin{eqnarray}\nonumber
\int_H \frac{1}{|x^0-\tilde{y}|^{n+\alpha}}dy&\geq&\int_D \frac{1}{|x^0-y|^{n+\alpha}}dy\\\nonumber
&=& \int_l^1\int_0^1\frac{\omega_{n-2}\tau^{n-2} d\tau}{(s^2+
\tau^2)^{\frac{n+\alpha}{2}}}ds\\\label{7}
&=& \int_l^1\int_0^{\frac{1}{s}}\frac{{\omega_{n-2}(st)^{n-2}}s dt}{s^{n
+\alpha}(1+t^2)^{\frac{n+\alpha}{2}}}ds\\\nonumber
&=& \int_l^1\frac{1}{s^{ 1+\alpha}}\int_0^{\frac{1}{s}}\frac{\omega_{n-2}t^{n-2}
dt}{(1+t^2)^{\frac{n+\alpha}{2}}}ds\\\nonumber
&\geq & \int_l^1\frac{1}{s^{1+\alpha}}\int_0^1\frac{\omega_{n-2}t^{n-2} dt}{(1+t^2)^{\frac{n+\alpha}{2}}}ds\\\label{8}
&\geq & C\int_l^1\frac{1}{s^{1+\alpha}}ds\rightarrow\infty,
\end{eqnarray}
where (\ref{7}) follows from the substitution $\tau=st$ and (\ref{8}) is
true when $l\rightarrow0$.

Hence $c(x)$ being lower bounded in $\Omega$ leads to
$$(-\lap)^{\alpha/2}u(x^0) +c(x^0)u(x^0)<0,\; \mbox{ when $l$ sufficiently small.}$$
This is a contradiction with condition (\ref{5}). Therefore, (\ref{6}) must be true.

\subsection{Decay at Infinity}

\begin{thm} Let $H=\{x \in \mathbb{R}^{n} | \; x_1<\lambda \, \mbox{ for some } \lambda \in \mathbb{R} \}$ and let $\Omega$ be an unbounded region in $H$.
Assume $u \in L_{\alpha} \cap C_{loc}^{1,1}(\Omega)$ is a solution of
\be
\left\{\begin{array}{ll}
(-\lap)^{\alpha/2}u(x) +c(x)u(x)\geq0  &\mbox{ in } \Omega,\\
u(x) \geq0&\mbox{ in } H\backslash\Omega,\\
u(\tilde{x})=-u(x)  &\mbox{ in } H,
\end{array}
\right.
\label{9}
\ee
with
\be
\underset{|x|\rightarrow \infty}{\underline{\lim}}|x|^{\alpha}c(x)\geq0,\label{9.1}
\ee
then there exists a constant $R_0>0$ ( depending on $c(x)$, but is independent of $w$ ) such that if
\be
u(x^0)=\underset{\Omega}{\min}\;u(x)<0,\label{10}
\ee
 then
\be
|x^0|\leq R_0.\label{11}
\ee
\label{thmDecay}
\end{thm}

$\mathbf{Proof.}$ \quad It follows from (\ref{9}) and (\ref{10}) that
\begin{eqnarray*}
(-\lap)^{\alpha/2}u(x^0)&=&C_{n,\alpha}PV\int_{\mathbb{R}^n}\frac{u(x^0)-u(y)}
{|x^0-y|^{n+\alpha}}dy\\
&=&C_{n,\alpha}PV\left\{\int_H \frac{u(x^0)-u(y)}
{|x^0-y|^{n+\alpha}}dy+\int_{\tilde{H}} \frac{u(x^0)-u(y)}
{|x^0-y|^{n+\alpha}}dy \right\}\\
&=&C_{n,\alpha}PV\left\{\int_H \frac{u(x^0)-u(y)}
{|x^0-y|^{n+\alpha}}dy+\int_H \frac{u(x^0)-u(\tilde{y})}
{|x^0-\tilde{y}|^{n+\alpha}}dy\right\}\\
&=&C_{n,\alpha}PV\left\{\int_H \frac{u(x^0)-u(y)}
{|x^0-y|^{n+\alpha}}dy+\int_H \frac{u(x^0)+u(y)}
{|x^0-\tilde{y}|^{n+\alpha}}dy\right\}\\
&\leq&
C_{n,\alpha}\int_H
\left\{\frac{u(x^0)-u(y)}{|x^0-\tilde{y}|^{n+\alpha}}
+\frac{u(x^0)+u(y)}{|x^0-\tilde{y}|^{n+\alpha}}\right\}dy\\
&=&C_{n,\alpha}\int_H \frac{2u(x^0)}
{|x^0-\tilde{y}|^{n+\alpha}}dy.
\end{eqnarray*}

For each fixed $\lambda$, when $|x^0| \geq \lambda$, we have
$B_{|x^0|}(x^1)\subset\tilde{H}$ with $x^1=(3|x^0|+x^0_1, (x^0)')$,  and it follows that
\begin{eqnarray*}
\int_H \frac{1}{|x^0-\tilde{y}|^{n+\alpha}}dy
&\geq&\int_{B_{|x^0|}(x^1)} \frac{1}{|x^0-y|^{n+\alpha}}dy\\
&\geq&\int_{B_{|x^0|}(x^1)} \frac{1}{4^{n+\alpha}|x^0|^{n+\alpha}}dy\\
&=&\frac{\omega_n}{4^{n+\alpha}|x^0|^{\alpha}}.
\end{eqnarray*}

Then we have
\begin{eqnarray*}
0&\leq &(-\lap)^{\alpha/2}u(x^0) +c(x^0)u(x^0)\\
&\leq&\left[\frac{2\omega_nC_{n,\alpha}}{4^{n+\alpha}|x^0|^{\alpha}}+
c(x^0)\right]u(x^0).
\end{eqnarray*}
Or equivalently,
$$\frac{2 \omega_nC_{n,\alpha}}{4^{n+\alpha} |x^0|^{\alpha}}+c(x^0)\leq 0.$$
Now if $|x^0|$ is sufficiently large, this would contradict (\ref{9.1}). Therefore, (\ref{11}) holds. This completes the proof.
\medskip

\begin{rem}

From the proof, one can see that the inequality
$$(-\lap)^{\alpha/2}u(x) +c(x)u(x)\geq0$$ and
condition (\ref{9.1}) are only required at points where $u$ is negative.

\end{rem}

\section{Method of Moving Planes and Its Applications}

\subsection{Radial Symmetry of $(-\lap)^{\alpha/2}u(x)=u^p(x),\; x \in R^n$}

\begin{thm}
Assume that $u \in   L_\alpha \cap C^{1,1}_{loc}$ and
\be
(-\lap)^{\alpha/2}u(x)=u^p(x),\; x \in R^n,\label{p1}
\ee
Then

(i) in the subcritical case $1<p<\frac{n+\alpha}{n-\alpha}$, (\ref{p1}) has no positive solution;

(ii) in the critical case $p=\frac{n+\alpha}{n-\alpha}$, the positive solutions must be radially symmetric and monotone decreasing about some point in $R^n$.
\end{thm}

\textbf{Proof.} Because no decay condition on $u$ near infinity is assumed, we are not able to
carry the {\em method of moving planes} on $u$ directly. To circumvent this difficulty, we make a Kelvin transform.

Let $x^0$ be a point in $\mathbb{R}^n$, and
$$ \bar{u}(x) = \frac{1}{|x-x^0|^{n-\alpha}} u\left( \frac{x-x^0}{|x-x^0|^2} + x^0 \right), \;\; x \in \mathbb{R}^n \setminus \{x^0\}. $$
be the Kelvin transform of $u$ centered at $x^0$. Then it is well-known that
\be
(-\lap)^{\alpha/2}\bar{u}(x)=\frac{\bar{u}^p(x)}{|x-x^0|^\tau}, \;\; x \in \mathbb{R}^n \setminus \{x^0\} \label{p3}
\ee
with $\tau=n+\alpha-p(n-\alpha)$. Obviously, $\tau = 0$ in the critical case.

Choose any direction to be the $x_1$ direction. For $\lambda < x^0_1$,
let $$T_\lambda=\{x \in R^n| \; x_1=\lambda\}, \quad x^\lambda=(2\lambda-x_1, x'),$$
$$ \bar{u}_\lambda(x)=\bar{u}(x^\lambda), \quad w_\lambda(x)=\bar{u}_\lambda(x)-\bar{u}(x), $$
and
$$\Sigma_\lambda=\{x \in R^n| x_1<\lambda\},\quad\tilde{\Sigma}_\lambda=\{x^\lambda|x \in \Sigma_\lambda\}.$$

First, notice that, by the definition of $w_{\lambda}$, we have
$$ \lim_{|x| \ra \infty} w_{\lambda}(x) = 0 .$$
Hence, if $w_{\lambda}$ is negative somewhere in $\Sigma_{\lambda}$,  then the negative minima of $w_{\lambda}$ were attained in the
interior of $\Sigma_{\lambda}$.

Let $$\Sigma_{\lambda}^- = \{ x \in \Sigma_{\lambda} \mid w_{\lambda}(x) < 0 \}.$$
Then from (\ref{p3}), we have, for $x \in \Sigma_{\lambda}^- \setminus \{(x^0)^{\lambda}\}$,
\begin{eqnarray*}
(-\lap)^{\alpha/2}w_\lambda(x)&=&\frac{\bar{u}_\lambda^p(x)}{|x^\lambda-x^0|^\tau}
-\frac{\bar{u}^p(x)}{|x-x^0|^\tau}\\
&\geq&\frac{\bar{u}_\lambda^p(x)-\bar{u}^p(x)}{|x-x^0|^\tau}\\
&\geq&\frac{p\bar{u}^{p-1}(x)w_\lambda(x)}{|x-x^0|^\tau};
\end{eqnarray*}
that is,
\be (-\lap)^{\alpha/2} w_\lambda(x) + c(x) w_\lambda(x) \geq 0 , \;\; x \in \Sigma_{\lambda}^- \setminus \{(x^0)^{\lambda}\},
\label{weq}
\ee
 with
 \be c(x) = - \frac{p\bar{u}^{p-1}(x)}{|x-x^0|^\tau}. \label{cx} \ee

\subsubsection{The Subcritical Case}
For $1<p<\frac{n+\alpha}{n-\alpha}$,
we show that (\ref{p1}) admits no positive solution.

\emph{Step 1.}

We show that, for $\lambda$ sufficiently negative,
\be
w_\lambda(x)\geq0 , \;\; \mbox{ in } \Sigma_\lambda \setminus \{(x^0)^{\lambda}\}.\label{p4}
\ee
This is done by using Theorem \ref{thmDecay} ({\em decay at infinity}).

First, we claim that for $\lambda$ sufficiently negative, there exists $\epsilon > 0$ and $c_{\lambda} >0$, such that
\begin{equation}
w_{\lambda} (x) \geq c_{\lambda}, \;\;\; \forall \, x \in B_{\epsilon} ((x^0)^{\lambda}) \setminus \{(x^0)^{\lambda}\} .
\label{bddaway}
\end{equation}
For the proof, please see the Appendix. From (\ref{bddaway}), one can see that $\Sigma_{\lambda}^-$ has no intersection with $B_{\epsilon}((x^0)^{\lambda})$

From (\ref{cx}), it is easy to verify that, for $|x|$ sufficiently large,
 \be c(x) \sim \frac{1}{|x|^{2\alpha}} . \label{cxsim} \ee
 Hence $c(x)$ satisfies condition (\ref{9.1}) in Theorem \ref{thmDecay}. Applying Theorem \ref{thmDecay} to $w_{\lambda}$ with
 $$H = \Sigma_{\lambda} \; \mbox{ and } \; \Omega =\Sigma_{\lambda}^-  $$
 for any sufficiently small $\epsilon$, we conclude that, there
 exists a $R_o >0$ (independent of $\lambda$), such that if
 $\bar{x}$ is a negative minimum of $w_{\lambda}$ in $\Sigma_{\lambda}$, then
 \be |\bar{x}| \leq R_o .
 \label{Ro}
 \ee

Now for $\lambda \leq -R_o,$  we must have
$$ w_{\lambda}(x) \geq 0 , \;\; \forall \, x \in \Sigma_{\lambda} \setminus \{(x^0)^{\lambda}\}.$$
This verifies (\ref{p4}).
\medskip

\emph{Step 2.} {\em Step 1} provides a starting point, from which we can now move the plane $T_{\lambda}$ to the right as long as (\ref{p4}) holds to its limiting position.

Let $$\lambda_0=\sup\{\lambda < x^0_1 \mid w_\mu(x)\geq0, \; \forall x \in \Sigma_\mu \setminus \{(x^0)^{\mu}\}, \mu\leq\lambda\}.$$
In this part, we show that $$\lambda_0=x^0_1$$ and
\be
 w_{\lambda_0}(x)\equiv0, \quad x \in \Sigma_{\lambda_0} \setminus \{(x^0)^{\lambda_0}\}.
 \ee

Suppose that $$\lambda_0<x^0_1, $$
we show that the plane $T_\lambda$ can be moved further right. To be more rigorous, there exists some
$\epsilon>0$, such that for any $\lambda \in (\lambda_0, \lambda_0+\epsilon)$, we have
\be
w_{\lambda}(x)\geq0, \quad x \in \Sigma_\lambda \setminus \{(x^0)^{\lambda}\}.\label{p10}
\ee
This is a contradiction with the definition of $\lambda_0$. Hence we must have
\be
\lambda_0 = x^0_1.  \label{p11}
\ee

Now we prove (\ref{p10}) by the combining use of  {\em narrow region principle} and  {\em decay at infinity}.
\medskip

 Again we need the fact (see the Appendix) that there exists $c_o > 0$ such that for sufficiently small $\eta$
 \begin{equation}
 w_{\lambda_0}(x) \geq c_o , \;\; \forall \, x \in B_{\eta}((x^0)^{\lambda_0}) \setminus \{(x^0)^{\lambda_0}\} .
 \label{bddaway0}
 \end{equation}

 By (\ref{Ro}), the negative minimum of $w_{\lambda}$ cannot be attained outside of $B_{R_o}(0)$. Next we argue that it can neither be attained inside of $B_{R_o}(0)$. Actually, we will show that for $\lambda$ sufficiently close to $\lambda_0$,
\be w_{\lambda} (x) \geq 0 , \;\; \forall \, x \in (\Sigma_{\lambda} \cap B_{R_o}(0)) \setminus \{(x^0)^{\lambda}\}.
\label{z7}
\ee

 From {\em narrow region principle} (Theorem \ref{thmNR}), there is a small $\delta >0$, such that for $\lambda \in [\lambda_0, \lambda_0 + \delta)$, if
\be
w_{\lambda} (x) \geq 0 , \;\; \forall \, x \in \Sigma_{\lambda_0 -\delta} \setminus \{(x^0)^{\lambda}\},
\label{z6}
\ee
then
\be
w_{\lambda} (x) \geq 0 , \;\; \forall \, x \in (\Sigma_{\lambda} \setminus \Sigma_{\lambda_0 -\delta}) \setminus \{(x^0)^{\lambda}\} .
\label{z1}
\ee
To see this, in Theorem \ref{thmNR}, we let
$$H = \Sigma_{\lambda}  \mbox{ and the narrow region } \Omega = (\Sigma_{\lambda}^- \setminus \Sigma_{\lambda_0 -\delta}),$$
while the lower bound of $c(x)$ can be seen from (\ref{cxsim}).
\medskip

Then what left is to show (\ref{z6}), and actually we only need
\be
w_{\lambda} (x) \geq 0 , \;\; \forall \, x \in (\Sigma_{\lambda_0 -\delta} \cap B_{R_o}(0)) \setminus \{(x^0)^{\lambda}\}.
\label{z8}
\ee

In fact, when $\lambda_0< x^0_1$, we have
\be
w_{\lambda_0}(x)>0, \quad x \in \Sigma_{\lambda_0} \setminus \{(x^0)^{\lambda_0}\}.\label{p12}
\ee
If not, there exists some $\hat{x}$ such that
$$w_{\lambda_0}(\hat{x})=\underset{\Sigma_{\lambda_0}}{\min}\;w_{{\lambda_0}}(x)=0.$$
It follows that
\begin{eqnarray}
(-\lap)^{\alpha/2}w_{\lambda_0}(\hat{x})
&=& C_{n,\alpha}PV\int_{\mathbb{R}^n}\frac{-w_{\lambda_0}(y)}{|\hat{x}-y|^{n+\alpha}}dy \nonumber\\
&=&C_{n,\alpha} PV \int_{\Sigma_{\lambda_0}} \frac{- w_{\lambda_0}(y)}{|\hat{x} - y|^{n+\alpha}}  dy + \int_{\mathbb{R}^n \setminus \Sigma_{\lambda_0}} \frac{- w_{\lambda_0}(y)}{|\hat{x} - y|^{n+\alpha}}  dy \nonumber \\
&=& C_{n,\alpha} PV \int_{\Sigma_{\lambda_0}} \frac{- w_{\lambda_0}(y)}{|\hat{x} - y|^{n+\alpha}}  dy + \int_{\Sigma_{\lambda_0}} \frac{w_{\lambda_0}(y)}{|\hat{x} - y^{\lambda}|^{n+\alpha}}  dy \nonumber \\
&=&  C_{n,\alpha} PV \int_{\Sigma_{\lambda_0}} \left( \frac{1}{|\hat{x} - y^{\lambda}|^{n+\alpha}}  -  \frac{1}{|\hat{x} - y|^{n+\alpha}} \right)
w_{\lambda_0}(y) dy \nonumber \\
&\leq& 0 .
\label{z3}
\end{eqnarray}

On the other hand
$$(-\lap)^{\alpha/2}w_{\lambda_0}(\hat{x})=\frac{\bar{u}_{\lambda_0}^p(\hat{x})}{|\hat{x}^{\lambda_0}-
x^0|^\tau}-\frac{\bar{u}^p(\hat{x})}{|\hat{x}-x^0|^\tau} = \frac{\bar{u}^p(\hat{x})}{|\hat{x}^{\lambda_0}-
x^0|^\tau}-\frac{\bar{u}^p(\hat{x})}{|\hat{x}-x^0|^\tau} > 0 .$$
A contradiction with (\ref{z3}). This proves (\ref{p12}).
It follows from (\ref{p12}) that there exists a constant $c_o > 0$, such that
$$w_{\lambda_0}(x)\geq c_o , \quad x \in \overline{\Sigma_{\lambda_0-\delta}\cap B_{R_o}(0)}.$$
Since $w_{\lambda}$ depends on $\lambda$ continuously, there exists $\epsilon > 0$ and $\epsilon < \delta$, such that for all $\lambda \in (\lambda_0, \lambda_0 + \epsilon)$,
we have
\be
w_{\lambda}(x)\geq 0, \quad x \in \overline{\Sigma_{\lambda_0-\delta}\cap B_{R_o}(0)}.
\label{z4}
\ee

Combining (\ref{z1}), (\ref{Ro}), and (\ref{z4}), we conclude that for all $\lambda \in (\lambda_0, \lambda_0 + \epsilon)$,
\be
w_{\lambda}(x)\geq 0, \quad x \in \Sigma_{\lambda} \setminus \{(x^0)^{\lambda}\}.
\label{z5}
\ee
This contradicts the definition of $\lambda_0$. Therefore, we must have
$$ \lambda_0 = x^0_1 \; \mbox{ and } w_{\lambda_0} \geq 0 \; \forall \, x \in \Sigma_{\lambda_0}. $$
Similarly, one can move the plane $T_\lambda$ from the $+\infty$ to the left and show that
\be
w_{\lambda_0} \leq 0 \; \forall \, x \in \Sigma_{\lambda_0}.
\ee

 Now we have shown that
 $$
 \lambda_0=x^0_1 \; \mbox{ and } \; w_{\lambda_0}(x)\equiv0, \quad x \in \Sigma_{\lambda_0}.
 $$
This completes {\em Step 2}.

So far, we have proved that $\bar{u}$ is symmetric about the plane $T_{x^0_1}$. Since the $x_1$ direction can be chosen arbitrarily, we have actually shown that $\bar{u}$ is radially symmetric
about $x^0$.

For any two points $X^i \in \mathbb{R}^n$, $i=$1, 2. Choose $x^0$ to be the midpoint: $ x^0 = \frac{X^1 + X^2}{2}$. Since $\bar{u}$ is radially symmetric about $x^0$, so is $u$, hence $u(X^1)=u(X^1)$.
This implies that $u$ is constant. A positive constant function does not satisfy (\ref{p1}). This proves the
nonexistence of positive solutions for (\ref{p1}) when $1<p<\frac{n+\alpha}{n-\alpha}$.

\subsubsection{The Critical Case}

Let $\bar{u}$ be the Kelvin transform of $u$ centered at the origin, then
\be
(-\lap)^{\alpha/2}\bar{u}(x)=\bar{u}^p(x).\label{p17}
\ee

We will show that either $\bar{u}$ is symmetric about the origin or $u$ is symmetric about some point.

We still use the notation as in the subcritical case. Step 1 is entirely the same as that in the subcritical case, that is, we can show that for $\lambda$ sufficiently negative,
$$ w_{\lambda}(x) \geq 0 , \;\; \forall \, x \in \Sigma_{\lambda} .$$

Let
$$\lambda_0=\sup\{\lambda \leq 0 |w_\mu(x)\geq0, \; \forall x \in \Sigma_\mu, \mu\leq\lambda\}.$$

\textbf{Case (i). } $\lambda_0<0$. Similar to the subcritical case, one can show that
$$ w_{\lambda_0}(x) \equiv 0 , \;\; \forall \, x \in \Sigma_{\lambda_0}. $$
It  follows that 0 is not a singular point of $\bar{u}$ and hence
$$u(x)=O(\frac{1}{|x|^{n-\alpha}})\; \mbox{ when }\;|x|\rightarrow\infty.$$
This enables us to apply the method of moving plane to $u$ directly and show that
$u$ is symmetric about some point in $R^n$.

\textbf{Case (ii). } $\lambda_0=0$. Then by moving the planes from near $x_1 = + \infty$, we derive that $\bar{u}$ is symmetric about the origin, and so does $u$.

In any case, $u$ is symmetric about some point in $\mathbb{R}^n$.

This completes the proof.

\subsection{A Dirichlet Problem on a Half Space}

We investigate a Dirichlet problem involving the fractional Laplacian on an upper half space
$$\mathbb{R}_+^n = \{ x = (x_1, \cdots, x_n) \mid x_n > 0 \}.$$

Consider
\begin{equation}
\left\{\begin{array}{ll}
(-\lap)^{\alpha/2} u = u^p (x), & x \in \mathbb{R}_+^n, \\
u(x) \equiv 0 , & x \not{\in} \mathbb{R}_+^n .
\end{array}
\right.
\label{e1h1}
\end{equation}

\begin{thm}
Assume that $0 < \alpha < 2$ and $u \in L_{\alpha} \cap C^{1,1}_{loc}$ is a nonnegative solution of
problem (\ref{e1h1}). Then
in the subcritical and critical case $1<p \leq \frac{n+\alpha}{n-\alpha}$, $u \equiv 0$.
\label{thm1h}
\end{thm}

To prove this theorem,
again we make a Kelvin transform. In order that $\mathbb{R}_+^n$ is invariant under the transform, we put the center $x^o$ on the boundary $\partial \mathbb{R}_+^n$.

Let
$$ v_{x^o} (x) = \frac{1}{|x-x^o|^{n-\alpha}} u\left( \frac{x-x^o}{|x-x^o|^2} + x^o \right). $$
be the Kelvin transform of $u$ centered at $x^o$. Then it is well-known that
\be
(-\lap)^{\alpha/2} v_{x^o}(x)=\frac{v_{x^o}^p(x)}{|x-x^o|^\tau}, \;\; x \in \mathbb{R}_+^n.
\label{k1}
\ee
with $\tau=n+\alpha-p(n-\alpha)$. Obviously, $\tau = 0$ in the critical case.
\medskip

{\em The main ideas are as follows. }
\smallskip

In the critical case $p = \frac{n+\alpha}{n-\alpha}$, we consider two possibilities.

 (i) \emph{There is a point $x^o\in \partial \mathbb{R}^n_+$, such that $v_{x^o}(x)$ is bounded near $x^o$}. In this situation, $u \leq \frac{C}{1+|x|^{n-\alpha}}$ has the needed asymptotic behavior near infinity, hence we move the planes in the direction of $x_n$-axis to show that the solution $u$ is monotone increasing in $x_n$.

 (ii) \emph{For all $x^o\in \partial \mathbb{R}^n_+$,  $v_{x^o}(x)$ are unbounded near $x^o$.} In this situation, we move the planes in $x_1, \cdots, x_{n-1}$ directions to show that, for every $x^o$,  $v_{x^o}$ is axially symmetric about the line that is parallel to $x_n$-axis and
 passing through $x^o$. This implies further that $u$ depends on $x_n$ only.

 In the subcritical case, we only need to work on $v_{x^o}(x)$; and similar to the above possibility (ii), we show that for every $x^o$,  $v_{x^o}$ is axially symmetric about the line that is parallel to $x_n$-axis and
 passing through $x^o$, which implies again that $u$ depends on $x_n$ only.

 In both cases, we will be able to derive contradictions.

 \subsubsection{The Critical Case}

 We consider two possibilities.

 (i) \emph{There is a point $x^o\in \partial \mathbb{R}^n_+$, such that $v_{x^o}(x)$ is bounded near $x^o$}. In this situation, from the symmetric expression
 $$ u (x) = \frac{1}{|x-x^o|^{n-\alpha}} v_{x^o} \left( \frac{x-x^o}{|x-x^o|^2} + x^o \right), $$
 we see immediately that
 \begin{equation}
 u(x) \sim \frac{1}{|x|^{n-\alpha}} , \;\; \mbox{ near infinity.}
 \label{usim}
 \end{equation}
 Consequently, by Theorem \ref{thm2.1}, we have
 $$\mbox{ either } u(x) > 0 \mbox{ or } u(x) \equiv 0 , \;\; \forall \, x \in \mathbb{R}^n_+ .$$
 Hence in the following, we may assume that $u >0$ in $\mathbb{R}^n_+.$

 Now we carry on the method of moving planes on the solution $u$ along $x_n$ direction.

 Let $$T_{\lambda} = \{ x \in \mathbb{R}^n \mid x_n = \lambda \} , \;\; \lambda > 0 ,$$
 and $$ \Sigma_{\lambda} = \{ x \in \mathbb{R}^n \mid 0< x_n < \lambda \} . $$
 Let $$ x^{\lambda} = (x_1, \cdots, x_{n-1}, 2\lambda - x_n)$$ be the reflection of $x$ about the plane $T_{\lambda}$.

 Denote $w_{\lambda} (x) = u(x^{\lambda}) - u(x)$, and
 $$ \Sigma_{\lambda}^- = \{ x \in \Sigma_{\lambda} \mid w_{\lambda}(x) < 0 \}.$$
 Then
 \be (-\lap)^{\alpha/2} w_\lambda(x) + c(x) w_\lambda(x) \geq 0 , \;\; x \in \Sigma_{\lambda}^- ,
\label{weq1}
\ee
 with
 \be c(x) = - p u^{p-1}(x).  \label{cx1} \ee

 From this and (\ref{usim}), we see that $c(x)$ is bounded from below in $\Sigma_{\lambda}^-$, and
  \be \lim_{|x| \ra \infty} w_{\lambda}(x) = 0 \mbox{ and } c(x) \sim \frac{1}{|x|^{2\alpha}} \mbox{ for $|x|$ large}.
  \label{wcx}
  \ee
 It follows that we can apply the {\em narrow region principle} to conclude that for $\lambda$ sufficiently small,
 \be
 w_{\lambda} (x) \geq 0 , \;\; \forall x \in \Sigma_{\lambda} ,
 \label{wxn}
 \ee
 because $\Sigma_{\lambda}$ is a narrow region.

 (\ref{wxn}) provides a starting point, from which we can move the plane $T_{\lambda}$ upward as long as inequality (\ref{wxn}) holds. Define
 $$ \lambda_o = \sup \{ \lambda \mid w_{\mu} (x) \geq 0, x \in \Sigma_{\mu}; \mu \leq \lambda \}.$$
 We show that
 \be \lambda_o = \infty .
 \label{y1}
 \ee

 Otherwise, if $\lambda_o < \infty$, then by (\ref{wcx}), combining the {\em Narrow Region Principle} and {\em Decay at Infinity} and going through the similar arguments as in the previous subsection, we are able to show that
 $$ w_{\lambda_o}(x) \equiv 0 \;\; \mbox{ in } \Sigma_{\lambda_o} ,$$
 which implies
 $$ u(x_1, \cdots, x_{n-1}, 2\lambda_o) = u(x_1, \cdots, x_{n-1}, 0) = 0 . $$
 This is impossible, because we assume that $u>0$ in $\mathbb{R}^n_+.$

 Therefore, (\ref{y1}) must be valid. Consequently, the solution $u(x)$ is monotone increasing with respect to $x_n$. This contradicts (\ref{usim}). Therefore what left to be considered is

 Possibility (ii): \emph{For all $x^o\in \partial \mathbb{R}^n_+$,  $v_{x^o}(x)$ are unbounded near $x^o$.}

 In this situation, we carry on the method of moving planes on $v_{x^o}$ along any direction in $\mathbb{R}^{n-1}$--the boundary of $\mathbb{R}^n_+$, call it $x_1$ direction.

For a given real number $\lambda$, define
$$ \hat{T}_{\lambda} = \{ x \in \mathbb{R}^n \mid x_1 = \lambda \} ,$$
$$
\hat{\Sigma}_{\lambda}=\{x=(x_{1},\cdots,x_{n})\in
\mathbb{R}^n_+ \mid x_{1}<\lambda\}
$$
and let
$$
x^{\lambda}=(2\lambda-x_{1},x_{2},\cdots,x_{n}).
$$

Let $w_{\lambda}(x) = v_{x^o}(x^{\lambda}) - v_{x^o}(x)$ and
$$ \hat{\Sigma}_{\lambda}^- = \{ x \in \hat{\Sigma}_{\lambda} \mid w_{\lambda}(x) < 0 \}.$$
 Then
 $$ (-\lap)^{\alpha/2} w_\lambda(x) + c(x) w_\lambda(x) \geq 0 , \;\; x \in \hat{\Sigma}_{\lambda}^- ,
$$
 with
 $$ c(x) = - p v_{x^o}^{p-1}(x). $$

 By the asymptotic behavior
 $$ v_{x^o}(x) \sim \frac{1}{|x|^{n-\alpha}} , \; \mbox{ for } |x|  \mbox{ large},$$
 we derive
 $$ \lim_{|x| \ra \infty} w_{\lambda} (x) = 0 \mbox{ and } c(x) \sim \frac{1}{|x|^{2\alpha}} \mbox{ for $|x|$ large}. $$
 These guarantee that we can apply the {\em narrow region principle} and {\em decay at infinity} to show the following:

 (i) For $\lambda$ sufficiently negative,
 $$ w_{\lambda} (x) \geq 0 , \;\; x \in \Sigma_{\lambda} .$$

(ii) Define
$$ \lambda_o = \sup \{ \lambda \mid w_{\mu} (x) \geq 0, x \in \Sigma_{\mu}; \mu \leq \lambda < x^o_1 \},$$
where $x^o_1$ is the first component of $x^o$. Then if $\lambda_o < x^o_1$, we must have
$$w_{\lambda_o}(x) \equiv 0 , \;\; x \in \Sigma_{\lambda_o} , $$
that is
$$ v_{x^o} (x^{\lambda_o}) \equiv v_{x^o}(x) , \;\; x \in \Sigma_{\lambda_o} .$$
This is impossible, because by our assumption, $v_{x^o}$ is unbounded near $x^o$, while it is bounded near $(x^o)^{\lambda_o}$. Therefore, we must have
$$\lambda_o = x^o_1 .$$

Based on this, and by moving the plane $\hat{T}_{\lambda}$ from near $x_1 = +\infty$ to the left to its limiting position, we show that $v_{x^o}$ is symmetric about the plane $\hat{T}_{x^o_1}$.
Since $x_1$ direction can be chosen arbitrarily, we conclude that $v_{x^o}$ is axially symmetric about the line parallel to $x_n$ axis and passing through $x^o$. Because $x^o$ is any point on $\partial \mathbb{R}^n_+$, we deduce that the original solution $u$ is independent of the first
$n-1$ variables, i.e, $u = u(x_n)$.

To finally derive a contradiction, we need two results from \cite{CFY}:

\begin{pro} (Theorem 4.1 in \cite{CFY})
Assume that $u \in L_{\alpha}$ is a locally bounded positive solution of
$$
\left\{\begin{array}{ll}(-\Delta)^{\alpha/2}u(x)=u^{p}(x),&x \in \mathbb{R}^n_+,\\
u(x)=0, & x \not\in \mathbb{R}_+^n .
\end{array}
\right.
$$
Then it is also a solution of
$$
u(x) = \int_{\mathbb{R}^n_+} G_{\infty} (x,y) u^p (y) dy ;
$$
and vice versa. Here $G_{\infty}(x,y)$ is the Green's function of the corresponding problem:
$$ G_{\infty}(x,y) = \frac{A_{n,\alpha}}{s^{\frac{n-\alpha}{2}}}\left[1-
B\frac{1}{(t+s)^{\frac{n-2}{2}}}
\int_{0}^{\frac{s}{t}}\frac{(s-tb)^{\frac{n-2}{2}}}{b^{\alpha/2}(1+b)}db\right],$$
with $s =|x-y|^2$ and $t = 4x_ny_n$.
\end{pro}

\begin{pro}
If $u=u(x_n) >0$, then
$$ \int_{\mathbb{R}^n_+} G_{\infty} (x,y) u^p (y) dy = \infty.$$
\end{pro}
(See the proof between page 23 and 27 in \cite{CFY}.)

Now these two propositions imply that if $u = u(x_n)$ is a positive solution of problem (\ref{e1h1}), then
$$ u(x) = \int_{\mathbb{R}^n_+} G_{\infty} (x,y) u^p (y) dy = \infty ,$$
which is obviously impossible.
This completes the proof in the critical case.

\subsubsection{The Subcritical Case}

 Recall that
 \be
(-\lap)^{\alpha/2} v_{x^o}(x)=\frac{v_{x^o}^p(x)}{|x-x^o|^\tau}, \;\; x \in \mathbb{R}_+^n.
\label{k2}
\ee

Similar to the possibility (ii) in the critical case, we carry on the method of moving planes on $v_{x^o}$ along any direction in $\mathbb{R}^{n-1}$, the boundary of $\mathbb{R}^n_+$, call it $x_1$ direction. Due to the presence of the term $\frac{1}{|x-x^o|^\tau}$ in equation
(\ref{k2}) with $\tau >0$, through a similar argument, we can derive that $v_{x^o}$ is axially symmetric about the line parallel to $x_n$ axis and passing through $x^o$, and hence the original solution $u$ is independent of the first $n-1$ variables, i.e, $u = u(x_n)$, which leads to a contradiction as in the critical case.

This completes the proof of the theorem.

\subsection{The Nonlinear Schr$\ddot{o}$diger Equation}
We study positive solutions for the nonlinear Schr$\ddot{o}$dinger equation with fractional diffusion
\begin{equation}
(-\lap)^{\alpha/2} u + u = u^p , \;\; x \in R^n .
\label{se1}
\end{equation}

\begin{thm}
Assume that $u \in L_{\alpha}\cap C_{loc}^{1,1}$ is a positive solution of (\ref{se1}) with $1<p<\infty$. If
\begin{equation}
\lim_{|x| \ra \infty} u(x) = a < \left( \frac{1}{p} \right)^{\frac{1}{p-1}} ,
\label{sec}
\end{equation}
then $u$ must be radially symmetric and monotone decreasing about some point in $\mathbb{R}^n$.
\label{thm3.2}
\end{thm}

{\bf Proof.} Because of the presence of the term $u$ in the equation, if one makes a Kelvin transform, the coefficients in the resulting equation do not possess the monotonicity needed in the method of moving planes. Hence we directly work on original equation (\ref{se1}).

Let $T_{\lambda}, \Sigma_{\lambda}, x^{\lambda}$, and $u_{\lambda}$  be defined as in the previous section. And let $w_{\lambda}(x) = u_{\lambda}(x) - u(x)$. Then at points where $w_{\lambda}$ is negative, it is easy to verify that
\begin{equation}
(-\lap)^{\alpha/2} w_{\lambda} + \left(1-pu^{p-1}\right) w_{\lambda} (x) \geq 0 .
\label{s1}
\end{equation}

{\em Step 1.} We apply Theorem \ref{thmDecay} ({\em decay at infinity})  to show that for sufficiently negative $\lambda$, it holds
\begin{equation}
w_{\lambda}(x) \geq 0 , \;\; x \in \Sigma_{\lambda} .
\label{w1}
\end{equation}
Here in (\ref{s1}), our $c(x) = \left(1-pu^{p-1}(x)\right)$.

First, by our assumption that $\lim_{|x| \ra \infty} u(x) = a$, we have, for each fixed $\lambda$,
$$ \lim_{|x| \ra \infty} w_{\lambda}(x) = 0 .$$
Hence if (\ref{w1}) is violated, then a negative minimum of $w_{\lambda}$ is attained at some
point, say at $x^o$.

By condition (\ref{sec}), we have,
$$ c(x) \geq 0 , \; \mbox{ for $|x|$ sufficiently large} , $$
and hence assumption (\ref{9.1}) in Theorem \ref{thmDecay} is satisfied. Consequently, there exists
$R_o$ (independent of $\lambda$), such that
\begin{equation}
 |x^o| \leq R_o .
 \label{xo}
 \end{equation}
It follows that, for $\lambda < -R_o$, we must have
$$ w_{\lambda}(x) \geq 0 , \;\; x \in \Sigma_{\lambda} .$$

{\em Step 2. } {\em Step 1} provides a starting point, from which we can now move the plane $T_{\lambda}$ to the right as long as (\ref{w1}) holds to its limiting position.

Let
$$\lambda_0=\sup\{\lambda \mid w_\mu(x)\geq0, \; \forall x \in \Sigma_\mu, \mu\leq\lambda\}.$$
It follows from (\ref{xo}) that $\lambda_0 < \infty .$

We will show that
\be
 w_{\lambda_0}(x)\equiv0, \quad x \in \Sigma_{\lambda_0}.
 \label{w=0}
 \ee

Suppose in the contrary,
$$ w_{\lambda_0}(x) \geq 0 \mbox{ and }  w_{\lambda_0}(x) \not{\equiv} 0 , \;  \mbox{ in } \Sigma_{\lambda_0}, $$
we must have
\begin{equation}
w_{\lambda_0}(x) >  0 \mbox{ in } \Sigma_{\lambda_0} .
\label{w>0}
\end{equation}
In fact, if (\ref{w>0}) is violated, then there exists a point $\hat{x} \in \Sigma_{\lambda_0}$, such that
$$ w_{\lambda_0}(\hat{x}) = \min w_{\lambda_0} = 0 .$$
Consequently, similar to (\ref{z3}), we have
\begin{eqnarray*}
(- \lap)^{\alpha/2} w_{\lambda_0}(\hat{x})
&=&  C_{n,\alpha} PV \int_{\Sigma_{\lambda_0}} \left( \frac{1}{|\hat{x} - y^{\lambda_0}|^{n+\alpha}}  -  \frac{1}{|\hat{x} - y|^{n+\alpha}} \right)
w_{\lambda_0}(y) dy \\
&<& 0 .
\end{eqnarray*}
This contradicts (\ref{s1}). Hence (\ref{w>0}) holds.

Then we show that the plane $T_\lambda$ can be moved further right. To be more rigorous, there exists some
$\epsilon>0$, such that for any $\lambda \in (\lambda_0, \lambda_0+\epsilon)$, we have
\be
w_{\lambda}(x)\geq0, \quad x \in \Sigma_\lambda.
\label{wgeq0}
\ee
This is a contradiction with the definition of $\lambda_0$. Therefore (\ref{w=0}) must be valid.

Under our assumptions, we have
$$ \lim_{|x| \ra \infty} w_{\lambda}(x) = 0 \mbox{ and $c(x)$ is bounded from below.} $$
Then combining the {\em narrow region principle} and the {\em decay at infinity}, through a similar argument as in the previous section, we derive (\ref{wgeq0}). This completes the proof of the theorem.

\subsection{More General Nonlinearities on a Bounded Domain}

Consider
\begin{equation}
\left\{ \begin{array}{ll}
(-\lap)^{\alpha/2} u(x)  = f(u(x)) , & x \in B_1(0), \\
u(x) = 0 , & x \not{\in} B_1(0).
\end{array}
\right.
\label{fu1}
\end{equation}
We prove
\begin{thm}
Assume that $u \in L_{\alpha}\cap C_{loc}^{1,1}(B_1(0))$ is a positive solution of (\ref{fu1}) with $f(\cdot)$ being Lipschitz continuous. Then $u$ must be radially symmetric and monotone decreasing about the origin.
\label{thmfu}
\end{thm}

{\bf Proof.}  Let $T_{\lambda}, x^{\lambda}, u_{\lambda},$ and $w_{\lambda}$ as defined in the previous section. Let
$$ \Sigma_{\lambda} = \{ x \in B_1(0) \mid x_1 < \lambda \} .$$
Then it is easy to verify that
$$ w_{\lambda}(x) + c_{\lambda}(x) w_{\lambda}(x) = 0 , \;\; x \in \Sigma_{\lambda}, $$
where
$$ c_{\lambda}(x) = \frac{f(u(x)) - f(u_{\lambda}(x))}{u(x) - u_{\lambda}(x)} .$$

Our Lipschitz continuity assumption on $f$ guarantees that $c_{\lambda}(x)$ is uniformly bounded
from below. Now we can apply Theorem \ref{thmNR} ({\em narrow region principle}) to conclude that
for $\lambda > -1$ and sufficiently close to $-1$,
\begin{equation}
w_{\lambda}(x) \geq 0 , \;\; x \in \Sigma_{\lambda} ;
\label{wgeq01}
\end{equation}
because $\Sigma_{\lambda}$ is a narrow region  for such $\lambda$.

Define
$$ \lambda_0 = \sup \{ \lambda \leq 0 \mid w_{\mu} (x) \geq 0, x \in \Sigma_{\mu}; \mu \leq \lambda \}.$$
Then we must have $\lambda_0 = 0$. Otherwise, we can use the {\em narrow region principle} and  similar arguments as in the previous section to show that we would be able to move the plane $T_{\lambda}$ further to the right to contradict the definition of $\lambda_0$. Therefore
$$
w_{0}(x) \geq 0 , \;\; x \in \Sigma_{0} ;
$$
or more apparently,
\begin{equation}
u(-x_1, x_2, \cdots, x_n) \leq u(x_1, x_2, \cdots, x_n) , \;\; 0<x_1<1.
\label{u0}
\end{equation}

Since the $x_1$-direction can be chosen arbitrarily, (\ref{u0}) implies $u$ is radially symmetric about the origin. The monotonicity is a consequence of the fact that (\ref{wgeq01}) holds for all
$-1<\lambda \leq 0$. This completes the proof of the theorem.

\section{Appendix}

Here, we prove (\ref{bddaway}) and (\ref{bddaway0}). Without loss of generality, we let $x^0 =0$.

\begin{lem}
Assume that $u \in C_{loc}^{1,1}( \mathbb{R}^n)$
 is a positive solution for
\begin{equation*}
  (-\lap)^{\alpha/2}u(x)=u^{p}(x), \quad x \in \mathbb{R}^n.
\end{equation*}
Let $v(x)=\frac{1}{|x|^{n-2}}u(\frac{x}{|x|^2})$ be the Kelvin
transform of $u$
 and
\begin{equation*}
  w_\la(x)=v(x^\la)-v(x).
\end{equation*}
Then there exists a constant $C>0$ such that for $\la$ sufficiently negative,
\begin{equation*}
 w_\la(x)\geq C>0, \quad x \in B_\varepsilon(0^\la)\backslash\{0^\la\}.
\end{equation*}
\end{lem}

Proof. From \cite{ZCCY}, we know that $u$ also satisfies the integral equation
$$ u(x) =  \int_{R^n}\frac{u^p(y)}{|x-y|^{n-\alpha}}dy  .$$

It follows that, for $|x|$ sufficiently large,
\begin{eqnarray*}
  u(x) &\geq& \int_{B_1(0)}\frac{u^p(y)}{|x-y|^{n-\alpha}}dy \\
   &\geq& \underset{B_1(0)}{\inf}\, u^p(x)\int_{B_1(0)}\frac{1}{|x-y|^{n-\alpha}}dy  \\
   &\geq & \frac{ c}{|x|^{n-\alpha}}.
\end{eqnarray*}
This implies that for small $\varepsilon>0$,
\begin{equation}\label{bdd7}
  v(x)\geq c, \quad x \in B_{\varepsilon}(0)\backslash\{0\}.
\end{equation}
When $\la$ is sufficiently negative, it holds that
\begin{equation}\label{bdd8}
 v(x)< \frac{c}{2}, \quad x \in B_\varepsilon(0^\la)\backslash\{0^\la\}.
\end{equation}
Combining (\ref{bdd7}) and (\ref{bdd8}), we arrive at
\begin{equation*}
 w_\la(x)\geq \frac{c}{2} > 0, \quad x \in B_\varepsilon(0^\la)\backslash\{0^\la\}.
\end{equation*}
This completes the proof of the Lemma.

\begin{lem} Let $w_{\la}$ be defined as before in the previous lemma and let
$$\lambda_0 = \sup\{ \la \mid w_{\mu}(x) \geq 0 , \forall \, x \in \Sigma_{\mu} \setminus \{0^{\mu}\}, \mu \leq \la \}.$$
If $w_{\lambda_0} \not{\equiv} 0$, then
\begin{equation*}
 w_{\la_0}(x)\geq c >0, \quad x \in B_\varepsilon(0^{\la_0})\backslash\{0^{\la_0}\}.
\end{equation*}
\end{lem}

Proof. By the definition, if $$w_{\la_0}(x)\not\equiv 0, \quad x \in \Sigma_{\la_0},$$
then there exists a point $x^0$ such that
$$w_{\la_0}(x^0)> 0.$$
And further, there exists a small positive $\delta$ such that
\begin{equation*}
 w_{\la_0}(x^0)\geq C_1> 0, \quad x \in B_\delta(x^0).
\end{equation*}
It follows that
\begin{eqnarray*}
w_{\la_0}(x)&=&\int_{\sum_{\la_0}} \bigg(\frac{1}{|x-y|^{n-\alpha}}-
  \frac{1}{|x-y^{\la_0}|^{n-\alpha} }\bigg)\big(v^\tau(y^{\la_0})-v^\tau(y) \big)dy\\
   &\geq& \int_{ B_\delta(x^0)} \bigg(\frac{1}{|x-y|^{n-\alpha}}-
   \frac{1}{|x-y^{\la_0}|^{n-\alpha} }\bigg)\big(v^\tau(y^{\la_0})-v^\tau(y) \big)dy\\
 &\geq& \int_{ B_\delta(x^0)} C_2\, C_1\,dy\\
  &\geq& C_3>0.
\end{eqnarray*}
This completes the proof.

\bigskip

{\em Authors' Addresses and E-mails:}
\medskip

Wenxiong Chen

Department of Mathematical Sciences

Yeshiva University

New York, NY, 10033 USA

wchen@yu.edu
\medskip

Congming Li

Department of Applied Mathematics

University of Colorado,

Boulder CO USA

cli@clorado.edu
\medskip

Yan Li

Department of Mathematical Sciences

Yeshiva University

New York, NY, 10033 USA

yali3@mail.yu.edu

\end{document}